\documentclass[12pt,singlespace,oneside]{article}
\usepackage{amsmath}
\usepackage{latexsym}

\input{tcilatex}
\begin{document}

\title{Helmholtz-Hodge Theorems: Unification of\ Integration and
Decomposition Perspectives\thanks{%
To Bill Grimm for his invaluable support, encouragement and wise advice over
the years.}}
\author{Jose G. Vargas \thanks{%
PST Associates, 138 Promontory Rd, Columbia, SC 29209-1244. USA,\newline
\hspace*{0.5cm} josegvargas@earthlink.net} \ }
\date{May 8, 2014.}
\maketitle

\begin{abstract}
We develop a Helmholtz-like theorem for differential forms in Euclidean
space $E_{n}$ using a uniqueness theorem similar to the one for vector
fields. We then apply it to Riemannian manifolds, $R_{n}$, which, by virtue
of the Schl\"{a}fli-Janet-Cartan theorem of embedding, are here considered
as hypersurfaces in $E_{N}$ with $N\geq n(n+1)/2$.

We obtain a Hodge decomposition theorem that includes and goes beyond the
original one, since it specifies the terms of the decomposition.

We then view the same issue from a perspective of integrability of the
system ($d\alpha =\mu ,$ $\delta \alpha =\nu $), relating boundary
conditions to solutions of ($d\alpha =0,$ $\delta \alpha =0$), [$\delta $ is
what goes by the names of divergence and co-derivative, inappropriate for
the Kaehler calculus with which we obtain the foregoing).
\end{abstract}

\section{Introduction}

We use the K\"{a}hler calculus (KC) \cite{K60}, \cite{K61}, \cite{K62} to
extend Helmholtz theorem to differential $k-$forms in Euclidean spaces of
arbitrary dimension. We then adapt it to Riemannian manifolds, which we view
as embedded in Euclidean spaces. We thus achieve a Hodge theorem which, not
only speaks of the decomposition of differential forms, but actually
specifies the terms of such decomposition as integrands involving the
exterior derivative and the co-differential of the differential form.

This paper has three major parts. In a first part we present the strategy to
get into the deep results we have announced. We then recollect formulas from
KC needed to get to those results. In a second part, we develop Helmholtz
theorem for differential $k-$forms in Euclidean spaces of arbitrary
dimension. In a third part, we adapt our most general Helmholtz theorem to
Riemannian manifolds, where the condition of vanishing differential $k-$form
at the boundary does not apply in general.

In the second part, the venue is a Euclidean space, $E_{n}$, and the
boundary condition is at infinity, where the differential form in question
is supposed to vanish If it does not, we assume its vanishing and add to the
result a constant differential needed for that vanishing.

In the third part, the venue is a region of a $E_{n}$, or a Riemannian
manifold; in any of the last two cases, dealing with the boundary condition
brings to the fore the emergence of a harmonic term which is additional to
the two terms characteristic of Helmholtz theorems. From now on and for
practical reasons, the appellative Helmholtz will be used whenever the
specifics of a problem lead to expressing a differential $k-$form as sum of
closed and co-closed terms exclusively. The appellative Hodge will be used
when a harmonic term also enters the expression.

It is important to emphasize two main issues which, if not made explicit,
may cause confusion. Helmholtz theorem is about integration. Hodge's theorem
is only about decomposition. In the first case, we speak of \textit{%
uniqueness of the solution }(of a differential system) that satisfies
certain conditions. In the second case, we speak of the \textit{uniqueness
of the decomposition} of a differential form into closed, co-closed and
harmonic terms. For differential $k-$forms and only them, we go beyond Hodge
by considering it from a perspective of integration and obtaining all three
terms in the decomposition as integrals. The decomposition theorem is, for
both differential $k-$forms and inhomogeneous ones, a by-product of the
integration theorem for differential $k-$forms.

The second issue has to do with the nature of differential forms in KC,
where components have three types of indices in the most general case. One
type is constituted by valuedness superscripts, as usual. One type of
subscripts is for components of multilinear functions of vector fields,
whether antisymmetric or not. Covariant differentiation governs
differentiation involving indices of those types. A second type of subscript
is for functions of $r-$surfaces, their evaluation being given by their
integration; exterior differentiation pertains to them. Thus, for example, $%
dx$ and $dx\mathbf{i}+dy\mathbf{j}+dz\mathbf{k}$ are functions of curves
whose evaluation on a curve $\gamma $ with end points $A$ and $B$ in $E_{3}$
is given by%
\begin{equation}
\int_{x_{A}}^{x_{B}}dx=x_{B}-x_{A},\text{ \ \ \ \ \ \ \ \ \ \ \ }%
\int_{A}^{B}d\mathbf{r}=(x_{B}-x_{A})\mathbf{i}+(y_{B}-y_{A})\mathbf{j}%
+(z_{B}-z_{A})\mathbf{k}.
\end{equation}%
Differential $r-$forms are here integrands (functions of $r-$surfaces) like
in Cartan \cite{Cartan23}, K\"{a}hler and Rudin \cite{Rudin}.

The view of differential forms just expressed is particularly suited for a
convenient and deep use of Laplacians for doing what Dirac delta functions
do in \textit{ad hoc} manner. Recall in this regard that the existence of a
Stokes theorem allows for the definition of their exterior derivatives, as
expressed by \'{E}. Cartan in 1922 (\cite{Cartan22}, section 74), thus years
before the concept of Dirac's delta. See also \cite{Cartan34}, sections 30
and 31.

Let $z$ represent the unit differential $n-$form in Euclidean space $E_{n}.$
Clearly $z=dr$ $r^{n-1}$ $\Omega _{n-1}$, where the differential $(n-1)-$%
form $\Omega _{n-1}$ is the unit element of ``solid angle'' and where $r$ is
the radial coordinate in $n-$dimensional space. Using a symbol $\partial $
to be later explained, we write the Laplacian as $\partial \partial .$ Given
a region $R$ (dimension $n$) of $E_{n}$ and containing the origin, one
applies Stokes theorem to the integral $(\partial \partial r^{-\lambda })z$
with undetermined integer $\lambda ,$ and where $r$ is the length of the
radius vector$.$ One thus has%
\begin{equation}
\int_{R}\left( \partial \partial \frac{1}{r^{\lambda }}\right)
z=-\int_{\partial R}\lambda r^{-(\lambda +1)}dr\,z=-\int_{\partial R}\lambda
\Omega _{n-1}\frac{r^{n-1}}{r^{\lambda +1}}
\end{equation}%
With $S_{n-1}$ as solid angle, the last term becomes $-(n-2)S_{n-1}$ if we
choose $\lambda $ to be equal to $n-2$. We take $-(n-2)S_{n-1}$ to the
opposite side in (2) to emphasize the independence of dimension:

\begin{equation}
1=-\text{ }\frac{1}{(n-2)S_{n-1}}\int_{R}\left( \partial \partial \frac{1}{%
r^{n-2}}\right) z.  \label{3}
\end{equation}

\section{Strategy}

K\"{a}hler based his calculus on K\"{a}hler algebra, i.e. the Clifford
algebra defined by%
\begin{equation}
dx^{i}dx^{j}+dx^{j}dx^{i}=2g^{ij}.  \label{4}
\end{equation}%
When at least one of two factors in a Clifford product is of grade one, it
can be decomposed as sum of exterior and interior products. Correspondingly,
K\"{a}hler's comprehensive derivative, here denoted as $\partial ,$ has two
pieces. We shall refer to them as exterior and interior derivatives ($d$ and 
$\delta $ respectively), given that they emerge through that composition of
the Clifford product (See equation (8)). Neither K\"{a}hler's comprehensive
derivative, nor its interior part, nor the exterior derivative of products
other than exterior ones satisfies the standard Leibniz rule (See next
section). Yet, as he did, we shall use the term derivative even in cases
when that rule is not satisfied.

We shall first provide the main formulas of the KC\ to be used in the paper.
We then proceed to do for differential $1-$forms in 3-D Euclidean space, $%
E_{3}$, what Helmholtz did for vector fields. Here the main difficulty to be
dealt with is that whereas the second exterior derivative is zero, the curl
of the curl of a vector field is not zero in general. So, there is not total
parallelism in the proofs pertaining to vector fields, on the one hand, and
differential $1-$forms, on the other. We have dealt with this issue in \cite%
{V53} and there is nothing to be added except for starting to use the
notation that will later be valid for arbitrary grade in $E_{n}.$

A Helmholz theorem for differential $2-$forms in $E_{3}$ is achieved by
expressing in the so obtained theorem the differential $1-$form in terms of
its dual. It is then simply a matter of solving for the differential $2-$%
form. Some needed polishing of the proof in \cite{V53} is given. Comparison
of those two theorems makes obvious their generalization to arbitrary grade
in arbitrary dimension.

Retrospectively, Helmholtz theorem will be seen as a particular case of our
better Hodge theorem, or the latter one as a generalization of the first.
The proofs of the intermediate versions of theorems in the evolution from
Helmholtz to Hodge are all cut by the same pattern. Because of the profuse
intermingling of appearances of the $d$ and $\delta $ operators, we proceed
to produce that general pattern.

We start, like Helmholtz did, with the hypothetical form of (in our case) a
differential form $\alpha $ as a sum of terms where those operators appear
in expressions of the form%
\begin{equation}
\alpha =d\left( ...\int ...\delta ^{\prime }\alpha ^{\prime }...\right)
+\delta \left( ...\int ...d^{\prime }\alpha ^{\prime }...\right) .
\end{equation}%
We are then led to compute only 
\begin{equation}
\delta d\left( ...\int ...\delta ^{\prime }\alpha ^{\prime }...\right) ,%
\text{ \ \ \ \ and \ \ \ \ }d\delta \left( ...\int ...d^{\prime }\alpha
^{\prime }...\right) .  \label{6}
\end{equation}%
since $dd$ and $\delta \delta $ vanish identically. The nature of the
contents of the integrands allows one to, following due process, introduce $%
\delta d$ and $d\delta $ inside the integral as $\delta ^{\prime }d^{\prime
} $ and $d^{\prime }\delta ^{\prime }.$ One then replaces $\delta ^{\prime
}d^{\prime }$ with $\partial ^{\prime }\partial ^{\prime }-d^{\prime }\delta
^{\prime }$, and $d^{\prime }\delta ^{\prime }$ with $\partial ^{\prime
}\partial ^{\prime }-\delta ^{\prime }d^{\prime }.$ The integration of the $%
\partial ^{\prime }\partial ^{\prime }$ terms are the ones that will produce 
$\delta \alpha $ and $d\alpha $ when applying $\delta $ and $d$ to the right
hand side of (5), respectively.

After doing that, we begin the transition to our richer Hodge theorem. It is
based on embedding Riemannian manifolds $R_{n}$, on Euclidean spaces, $%
E_{N}. $ Recall that the standard Hodge theorem has three terms, one of them
closed, another one co-closed and a third one which is harmonic, the last
one being zero in Helmholtz theorem. But it would not be zero if, instead of
integrating to the whole of $E_{n}$, and of $E_{3}$ in particular, we
integrated to regions thereof.

The next step in the argument consists in embedding a Riemannian manifold in
a Euclidean space. Such embedding gives us a $n-$hypersurface in a $N-$%
Euclidean space. a hypersurface is not a region, which has the same
dimension as the Euclidean space as a matter of definition of the term. By
the time we shall have reached this step, it will be obvious that this is
not a problem whatsoever, which would be if dealing with vector fields (we
shall discuss this in a later section). In the case of a region as in the
case of a hypersurface, it is the integration of the term containing $%
d^{\prime }\delta ^{\prime }$ (arising from the replacement of $\delta
^{\prime }d^{\prime }$ with $\partial ^{\prime }\partial ^{\prime }-$ $%
d^{\prime }\delta ^{\prime }$), and of the term $\delta ^{\prime }d^{\prime
} $ (arising from the replacement of $d^{\prime }\delta ^{\prime }$ with $%
\partial ^{\prime }\partial ^{\prime }-$ $\delta ^{\prime }d^{\prime }$)
that gives rise to the harmonic term in Hodge's theorem. The harmonic term
is an effect of the boundary condition. Said better, it is a reflection of
what happens outside the boundary.

\section{Basic\textbf{\ K\"{a}hler calculus}}

\subsection{\textbf{Concepts}}

For dealing with Riemannian spaces, we shall embed them in Euclidean spaces.
We may then use Cartesian coordinates to simplify the computations. We shall
use Roman characters for differential forms in equations from K\"{a}hler's
papers and Greek letters otherwise.

K\"{a}hler defines covariant derivatives of differential forms. In terms of
Cartesian coordinates, take the very simple form%
\begin{equation}
d_{h}v=\frac{\partial v}{\partial x^{h}}.
\end{equation}%
when they are scalar-valued. $d_{h}$ satisfies the Leibniz rule, 
\begin{equation}
d_{h}(u\vee v)=d_{h}u\vee v+u\vee d_{h}v,
\end{equation}%
As mentioned in the previous section,

\begin{equation}
\partial v\equiv dx^{h}\vee d_{h}v=dv+\delta v.
\end{equation}%
\begin{equation}
dv\equiv dx^{h}\wedge d_{h}v,\text{ \ \ \ \ }\delta v\equiv dx^{h}\cdot
d_{h}v.  \label{10}
\end{equation}%
$dv$ is the exterior derivative. For $\delta v$, see below.

It follows from (9) that, in principle, $\partial \partial $ consists of
four terms. One of them is $dd$, which, as we know, is zero. Similarly, $%
\delta \delta $ equals $0$ if $d_{h}$ is computed with the Levi-Civita
connection. This is automatically the case in Euclidean spaces. We thus have
the well known equation%
\begin{equation}
\partial \partial =\delta d+d\delta .
\end{equation}

\subsection{Differentiation of products}

Let operators $\eta $ and $e^{h}$ be distributive operators defined on $k-$%
forms as%
\begin{equation}
\eta u_{r}=(-1)^{r}u_{r}\text{, \ \ \ \ }e^{h}u=dx^{h}\cdot u,
\end{equation}%
The following two sets of equations%
\begin{equation}
\partial (u\vee v)=\partial u\vee v+\eta u\vee \partial v+2e^{h}u\vee d_{h}v,
\end{equation}%
\begin{equation}
d(u\vee v)=du\vee v+\eta u\vee dv+e^{h}u\vee d_{h}v-\eta d_{h}u\vee e_{h}v,
\end{equation}%
\begin{equation}
\delta (u\vee v)=\delta u\vee v+\eta u\vee \delta v+e^{h}u\vee d_{h}v+\eta
d_{h}u\vee e_{h}v
\end{equation}%
and%
\begin{equation}
\partial (u\wedge v)=\partial u\wedge v+\eta u\wedge \partial v+e^{h}u\wedge
d_{h}v+\eta d_{h}u\wedge e_{h}v,  \label{16}
\end{equation}%
\begin{equation}
d(u\wedge v)=du\wedge v+\eta u\wedge dv,\text{ \ \ \ \ \ \ \ \ \ \ \ \ \ \ \
\ \ \ \ \ \ \ \ \ \ \ \ \ \ \ \ \ }
\end{equation}%
\begin{equation}
\delta (u\wedge v)=\delta u\wedge v+\eta u\wedge \delta v+e^{h}u\wedge
d_{h}v+\eta d_{h}u\wedge e_{h}v
\end{equation}%
provide some of the flavor of differentiation in KC.

When using (8) to get Eq. (13), one has to pass $dx^{h}$ to the right of $u$
to multiply $d_{h}v.$ It is in this process that the last term arises. In
Eqs. (14)-(16) and (18), there are other terms arising from the interplay of
the products $\vee $ and $\wedge $, one of them explicit and the other one
implicit in each of $d(u\vee v)$ and $\partial (u\wedge v)$. Equations (15)
and (8) are the respective differences between (13) and (14), on the one
hand, and (16) and (17) on the other.

Of great importance is the concept of \textit{constant differential, }$c$,
defined by $d_{h}c=0.$ Then%
\begin{equation}
\partial (u\vee c)=(\partial u)\vee c,
\end{equation}%
but $\partial (c\vee u)\neq c\vee (\partial u).$ For differential $0-$forms, 
$f$, we have, using (18),%
\begin{equation}
\delta (cf)=(\eta c)\delta f+(dx^{h}\cdot c)d_{h}f=0+df\cdot c=-\eta c\cdot
df,
\end{equation}%
since $\delta f$ is zero and $d_{h}f$ is a $0-$form. Polynomials in
Cartesian $dx$'s with constant coefficients are constant differentials.

\subsection{About the interior derivative}

When the connection is Levi-Civita's, $\delta $ is the co-derivative. In $%
E_{n}:$\textbf{\ }%
\begin{equation}
\delta v=(-1)^{n(n-1)/2}zd(vz)\text{, \ \ \ \ \ \ \ \ }z\delta v=d(vz)\text{,%
}
\end{equation}%
where $z$ is the unit $n-$form for given orientation. Its square is $%
(-1)^{n(n-1)/2}$. Define $u\equiv $ $vz.$ Then $v$ equals $%
(-1)^{n(n-1)/2}uz. $ Multiply the second of (21) by $z.$ We get: 
\begin{equation}
zz\delta \lbrack (-1)^{n(n-1)/2}uz]=\delta (uz)=zdu.
\end{equation}%
We have obtained the following useful formulas%
\begin{equation}
zdu=\delta (uz)\text{, \ \ \ \ \ \ \ \ }z\delta u=d(uz)\text{,}  \label{23}
\end{equation}%
the last one being the same as in (21), after a change of notation.

\section{Helmholtz Uniqueness}

We now show that, under the usual conditions for Helmholtz theorem, the
solution is unique for differential $k-$forms, i.e. of defined grade (here
named homogeneous). The proof is based on K\"{a}hler's very comprehensive
version for differential forms \cite{K60}, \cite{K62} of Green's symmetric
theorem, also called second Green identity.

\subsection{K\"{a}hler's Green theorem}

Any differential form can be written as a sum of monomials. Assume each
monomial written as a product of differential 1-forms. The operator $\zeta $
will denote the reversion of all factors in such products.

Let the subscript zero denote the $0-$form part. K\"{a}hler defines the
scalar product of order zero as the differential $n-$form%
\begin{equation}
(u,v)\equiv (\zeta u\vee v)_{0}\text{ }wz=(\zeta u\vee v)\wedge w.
\end{equation}%
If $u$ and $v$ are homogeneous, it is necessary but not sufficient condition
that $r=s$ in $(u_{r},v_{s})$ to be different from zero. He defines the
scalar product as the differential $(n-1)-$form $(u,v)_{1}$ defined by%
\begin{equation}
(u,v)_{1}\equiv dx^{i}\cdot (dx^{i}\vee u,v).  \label{25}
\end{equation}%
He then proves the ``Green-K\"{a}hler theorem'':%
\begin{equation}
d(u,v)_{1}=(u,\partial v)+(v,\partial u)
\end{equation}

\subsection{Helmholtz uniqueness for differential $k-$forms}

Let ($u_{1}$, $u_{2})$ be differential $k-$forms such that $du_{1}=du_{2}$, $%
\delta u_{1}=\delta u_{2}$ on a differential manifold, $R$, and such that,
at the boundary, $\partial R,$ $u_{1}$ equals $u_{2}.$ Define $\beta
=u_{1}-u_{2}.$ Hence%
\begin{equation}
d\beta =0=\delta \beta \text{ on }R,\text{ \ \ \ \ \ \ \ \ \ \ \ \ \ \ \ \ }%
\beta =0\text{ on }\partial R,
\end{equation}%
and, locally,%
\begin{equation}
(\beta =d\alpha ,\text{ }\delta d\alpha =0)\text{ \ \ on }R,\text{ \ \ \ \ \
\ \ \ \ }d\alpha =0\text{ \ \ on }\partial R.
\end{equation}%
Equation (26) with $u=\alpha $ and $v=d\alpha $ reads%
\begin{equation}
d(\alpha ,d\alpha )_{1}=(\alpha ,\partial d\alpha )+(d\alpha ,\partial
\alpha ).
\end{equation}%
By (9) and (28), we have%
\begin{equation}
(\alpha ,\partial d\alpha )=(\alpha ,dd\alpha )+(\alpha ,\delta d\alpha
)=0+0.  \label{30}
\end{equation}%
Consider next $(d\alpha ,\partial \alpha )$ and use again (9). If $\alpha $
is of definite grade, so are $d\alpha $ and $\delta \alpha $, but their
grades differ by two units. Their scalar product is, therefore, zero. On the
other hand, we have, with $a_{A}$ defined by $d\alpha =a_{A}dx^{A}$ (with
summation over the algebra as a module),%
\begin{equation}
(d\alpha ,\partial \alpha )=(d\alpha ,\delta \alpha )+(d\alpha ,d\alpha
)=0+\sum \left| a_{A}\right| ^{2}.
\end{equation}%
Substituting (30)-(31) in (29), applying Stokes theorem and using $d\alpha
=0 $, we get 
\begin{equation}
0=\int_{\partial R}(\alpha ,d\alpha )_{1}=\int_{R}d(\alpha ,d\alpha
)_{1}=\int_{R}\sum \left| a_{A}\right| ^{2}.
\end{equation}%
Thus $0=d\alpha =\beta =u_{1}-u_{2}$, and, therefore, $u_{1}=u_{2}.$
Uniqueness under conditions like those for Helmholz theorem has been proved.

Let us not overlook that this theorem has been derived under the assumption
that $\alpha $ and, therefore, $\beta $ and the $u$'s are of definite
grades. It appears that there is no Helmholtz uniqueness theorem for
inhomogeneous differential forms. The equation $(d\alpha ,\delta \alpha )=0$
constitutes a formidable system of equations, likely to have infinite
solutions (even if the space were two dimensional!) This is nevertheless no
impediment to prove Hodge's standard decomposition theorem.

\section{Helmholtz Theorems for $k-$forms}

\subsection{Helmholtz Theorem for $k-$forms in $E_{3}$}

In this section, we shall try to avoid potential confusion by replacing the
symbol $z$ with the symbol $w$.

With $r_{12}\equiv \lbrack (x-x^{\prime })^{2}+(y-y^{\prime
})^{2}+(z-z^{\prime })^{2}]^{1/2}$, the standard Helmholtz theorem states%
\begin{equation}
\mathbf{v}=-\frac{1}{4\pi }\mathbf{\nabla }\int_{E_{3}^{\prime }}\frac{%
\mathbf{\nabla }^{\prime }\cdot \mathbf{v}(\mathbf{r}^{\prime })}{r_{12}}%
dV^{\prime }+\frac{1}{4\pi }\mathbf{\nabla }\times \int_{E_{3}^{\prime }}%
\frac{\mathbf{\nabla }^{\prime }\times \mathbf{v}(\mathbf{r}^{\prime })}{%
r_{12}}dV^{\prime }.
\end{equation}%
By proceeding in parallel to the proof of (33), we showed in \cite{V53} that
Helmholtz theorem for differential $1-$forms in $E_{3}$ reads%
\begin{equation}
\alpha =-\frac{1}{4\pi }d\int_{E_{3}^{\prime }}\frac{(\delta ^{\,\prime
}\alpha ^{\prime })w^{\prime }}{r_{12}}-\frac{1}{4\pi }\delta \left(
dx^{j}dx^{k}\int_{E_{3}^{\prime }}\frac{d^{\prime }\alpha ^{\prime }\wedge
dx^{\prime i}}{r_{12}}\right) ,
\end{equation}%
This theorem is a particular case of the theorem in the next\ subsection.
Needless to say that inputs $\delta \alpha $ and $d\alpha $ must be $0-$form
and $2-$forms.

The equation $\alpha =$ $w\beta $, uniquely defines $\beta .$ We substitute
it in (34) and solve for $\beta $:%
\begin{equation}
\beta =\frac{1}{4\pi }wd\left( \int_{E_{3}^{\prime }}\frac{\delta ^{\,\prime
}(w^{\prime }\beta ^{\prime })}{r_{12}}w^{\prime }\right) +\text{ }\frac{1}{%
4\pi }w\delta \left( dx^{jk}\int_{E_{3}^{\prime }}\frac{d^{\prime
}(w^{\prime }\beta ^{\prime })\wedge dx^{\prime i}}{r_{12}}\right) .
\label{35}
\end{equation}%
We use (18) to move $w$ to the right of $d$ and $\delta .$ For details, we
refer to \cite{V53}, except for the following simplifications.

Denote the first integral in (35) as $I$ and the second one as $I^{i}.$ We
have $wdI=\delta (wI)$ and $w^{\prime }\delta ^{\,\prime }(w^{\prime }\beta
^{\prime })=w^{\prime }(w^{\prime }d\beta ^{\prime })=-d\beta ^{\prime
}=-d\beta ^{\prime }\wedge 1.$ The exterior product by $1$ is superfluous,
except for the purpose for making later Eq. (39) clear. Similarly, $w\delta
(dx^{jk}I^{i})=d(wdx^{jk}I^{i})=-d(dx^{i}I^{i})$ and%
\begin{equation*}
d^{\prime }(w^{\prime }\beta ^{\prime })\wedge dx^{\prime i}=dx^{\prime
i}\wedge d^{\prime }(w^{\prime }\beta ^{\prime })=\frac{1}{2}\left[
dx^{\prime i}w^{\prime }\delta ^{\,\prime }\beta ^{\prime }+w^{\prime
}\delta ^{\,\prime }\beta ^{\prime }dx^{\prime i}\right] =
\end{equation*}%
\begin{equation}
\text{ \ \ \ \ \ \ \ \ \ \ }=\frac{1}{2}(dx^{\prime jk}\delta ^{\,\prime
}\beta ^{\prime }+\delta ^{\,\prime }\beta ^{\prime }dx^{\prime
jk})=dx^{\prime jk}\wedge \delta ^{\,\prime }\beta ^{\prime }=\delta
^{\,\prime }\beta ^{\prime }\wedge dx^{\prime jk}.
\end{equation}%
We use these results in (35), change the order of the terms and get%
\begin{equation}
\beta =-\text{ }\frac{1}{4\pi }d\left( dx^{i}\int_{E_{3}^{\prime }}\frac{%
\delta ^{\,\prime }\beta ^{\prime }\wedge dx^{\prime jk}}{r_{12}}\right) -%
\frac{1}{4\pi }\delta \left( w\int_{E_{3}^{\prime }}\frac{\delta ^{\,\prime
}\beta ^{\prime }\wedge 1}{r_{12}}\right) ,
\end{equation}

Write the first term in (34) as%
\begin{equation}
-\frac{1}{4\pi }d\left[ 1\wedge \int_{E_{3}^{\prime }}\frac{(\delta
^{\,\prime }\alpha ^{\prime })\wedge w^{\prime }}{r_{12}}\right] .
\end{equation}%
Let the index $A$ label a Cartesian basis of the algebra as module. Let $dx^{%
\bar{A}}$ be the unique element in the basis such that $dx^{A}\wedge dx^{%
\bar{A}}=w.$ Define $\int_{E_{3}}\gamma _{r}$ if the grade $r$ of $\gamma $
is different from $3.$ All four terms on the right of (34) and (37) are thus
of the form%
\begin{equation}
-\frac{1}{4\pi }d\left[ dx^{A}\int_{E_{3}^{\prime }}\frac{(\delta ^{\,\prime
}\_\_)\wedge dx^{\prime \bar{A}}}{r_{12}}\right] \text{ \ \ \ \ or \ \ \ \ }-%
\frac{1}{4\pi }\delta \left[ dx^{A}\int_{E_{3}^{\prime }}\frac{(d^{\prime
}\_\_)\wedge dx^{\prime \bar{A}}}{r_{12}}\right] .
\end{equation}%
Take, for instance, the first of the two expressions in (39). We sum over
all $A$, equivalently, over all $\bar{A}^{\prime }.$ The grade of $(\delta
^{\,\prime }\_\_)$ determines the grade of the only $dx^{\prime \bar{A}}$
that may yield not zero integral since the sum of the respective grades must
be $3.$ For each surviving value of the index $\bar{A}$, the value of the
index $A$ ---thus the specific $dx^{A}$ at the front of the integral--- is
determined. We shall later show for ulterior generalization that we may
replace the Cartesian basis with any other basis, which we shall choose to
be orthonormal since they are the ``canonical ones'' of Riemannian spaces.

\subsection{Helmholtz Theorem for Differential k-forms in $E_{n}$}

Let $\omega ^{A}$ ($\equiv \omega ^{i_{1}}\omega ^{i_{2}}...\omega ^{i_{r}}$%
) denote elements of a basis in the K\"{a}hler algebra of differential forms
such that the $\omega ^{\mu }$ are orthonormal. The purpose of using an
orthonormal basis is that exterior products can be replaced with Clifford
products. Let $\omega ^{\bar{A}}$ be the monomial (uniquely) defined by $%
\omega ^{A}\omega ^{\bar{A}}=z$, with no sum over repeated indices.

The generalized Helmholtz theorem in $E_{n}$ reads as follows%
\begin{equation}
\alpha =-\frac{1}{(n-2)S_{n-1}}[d(\omega ^{A}I_{A}^{\delta })+\delta (\omega
^{A}I_{A}^{d})],  \label{40}
\end{equation}%
with summation over a basis in the algebra and where%
\begin{equation}
I_{A}^{\delta }\equiv \int_{E_{n}^{\prime }}\frac{(\delta ^{\prime }\alpha
^{\prime })\wedge \omega ^{\prime \bar{A}}}{r_{12}^{n-2}},\text{ \ \ \ \ }%
I_{A}^{d}\equiv \int_{E_{n}^{\prime }}\frac{(d^{\prime }\alpha ^{\prime
})\wedge \omega ^{\prime \bar{A}}}{r_{12}^{n-2}}.
\end{equation}%
$r_{12}$ is defined by $r_{12}^{2}=(x_{1}-x_{1}^{\prime
})^{2}+...+(x_{n}-x_{n}^{\prime })^{2}$ in terms of Cartesian coordinates.

It proves convenient for performing differentiations to replace $\omega ^{i}$%
, $\omega ^{A}$ and $\omega ^{\bar{A}}$ with $dx^{i}$, $dx^{A}$ and $dx^{%
\bar{A}}$. If the results obtained are invariants, one can re-express the
results in terms of arbitrary bases.

We proceed again via the uniqueness theorem, as in the vector calculus, with
specification now of $d\alpha ,$ $\delta \alpha $ and that $\alpha $ goes
sufficiently fast at $\infty .$ vanishing of $\alpha $ at infinity. Because
of the annulment of $dd$ and $\delta \delta $, the proof reduces to showing
that $\delta d(dx^{A}I_{A}^{\delta })$ and $d\delta (dx^{A}I_{A}^{d})$
respectively yield $\delta \alpha $ and $d\alpha $, up to the factor at the
front in (40). Since the treatment of both terms is the same, we shall carry
them in parallel, as in%
\begin{equation}
\left( 
\begin{array}{c}
{\small \delta } \\ 
{\small d}%
\end{array}%
\right) \alpha \rightarrow \left( 
\begin{array}{c}
{\small \delta d} \\ 
{\small d\delta }%
\end{array}%
\right) dx^{A}I_{A}^{\tbinom{\delta }{d}}=\partial \partial dx^{A}I_{A}^{%
\tbinom{\delta }{d}}-\left( 
\begin{array}{c}
{\small d\delta } \\ 
{\small \delta d}%
\end{array}%
\right) dx^{A}I_{A}^{\tbinom{\delta }{d}}.
\end{equation}%
In the first term on the right hand side of (42), we move $\partial \partial 
$ to the right of $dx^{A}$, insert it inside the integral with primed
variables, multiply by $-\frac{1}{(n-2)S_{n-1}}$ and treat the integrand as
a distribution. We easily obtain that the first term yields $\tbinom{\delta
\alpha }{d\alpha }$.

For the last term in (42), we have%
\begin{equation}
\left( 
\begin{array}{c}
{\small d\delta } \\ 
{\small \delta d}%
\end{array}%
\right) dx^{A}I_{A}^{\tbinom{\delta }{d}}=\left( 
\begin{array}{c}
d\left[ dx^{i}\cdot dx^{A}\frac{\partial I_{A}^{\delta }}{\partial x^{i}}%
\right] \\ 
\delta \left[ (\eta dx^{A})\wedge dx^{i}\frac{\partial I_{A}^{d}}{\partial
x^{i}}\right]%
\end{array}%
\right) .
\end{equation}%
For the first line in (43), we have used (7) and the second equation (10).
For the development of the second line, we have used the Leibniz rule.

We use the same rule to also transform the first line in (43),%
\begin{equation}
d\left( dx^{i}\cdot dx^{A}\frac{\partial I_{A}^{\delta }}{\partial x^{i}}%
\right) =[\eta (dx^{i}\cdot dx^{A})]\wedge dx^{l}\frac{\partial
^{2}I_{A}^{\delta }}{\partial x^{l}\partial x^{i}}=(dx^{A}\cdot
dx^{i})\wedge dx^{l}\frac{\partial ^{2}I_{A}^{\delta }}{\partial
x^{l}\partial x^{i}}.
\end{equation}%
For the second line, we get%
\begin{equation}
\delta \left[ (\eta dx^{A})\wedge dx^{i}\frac{\partial I_{A}^{d}}{\partial
x^{i}}\right] =dx^{l}\cdot \left[ \frac{\partial ^{2}I_{A}^{d}}{\partial
x^{l}\partial x^{i}}(\eta dx^{A})\wedge dx^{i}\right] .  \label{45}
\end{equation}%
We shall use here that 
\begin{equation}
dx^{l}[(\eta dx^{A})\wedge dx^{i}]=-\eta \lbrack \eta (dx^{A}\wedge
dx^{i})]\cdot dx^{l}=(dx^{A}\wedge dx^{i})\cdot dx^{l},
\end{equation}%
thus obtaining%
\begin{equation}
\delta \left[ (\eta dx^{A})\wedge dx^{i}\frac{\partial I_{A}^{d}}{\partial
x^{i}}\right] =(dx^{A}\wedge dx^{i})\cdot dx^{l}\frac{\partial ^{2}I_{A}^{d}%
}{\partial x^{l}\partial x^{i}}.
\end{equation}%
Getting (44) and (47) into (43), we obtain%
\begin{equation}
\left( 
\begin{array}{c}
{\small d\delta } \\ 
{\small \delta d}%
\end{array}%
\right) dx^{A}I_{A}^{\tbinom{\delta }{d}}=\left[ dx^{A}{\LARGE (}_{\text{ }%
\wedge }^{\text{ }\cdot }{\Large )}dx^{i}\right] {\LARGE (}_{\text{ }\cdot
}^{\text{ }\wedge }{\LARGE )}dx^{l}\int_{E_{n}^{\prime }}\frac{\partial ^{2}%
}{\partial x^{\prime i}\partial x^{\prime l}}\frac{1}{r_{12}^{n-2}}\left( 
\begin{array}{c}
{\small \delta }^{\prime }{\small \alpha }^{\prime } \\ 
{\small d}^{\prime }{\small \alpha }^{\prime }%
\end{array}%
\right) \wedge dx^{\prime \bar{A}}.
\end{equation}

Integration by parts with respect to $x^{\prime i}$ yields two terms. The
total differential term is%
\begin{equation}
\left[ dx^{A}{\LARGE (}_{\text{ }\wedge }^{\text{ }\cdot }{\Large )}dx^{i}%
\right] {\LARGE (}_{\text{ }\cdot }^{\text{ }\wedge }{\LARGE )}%
dx^{l}\int_{E_{n}^{\prime }}\frac{\partial }{\partial x^{\prime i}}\left[
\left( \frac{\partial }{\partial x^{\prime l}}\frac{1}{r_{12}^{n-2}}\right)
\left( 
\begin{array}{c}
{\small \delta }^{\prime }{\small \alpha }^{\prime } \\ 
{\small d}^{\prime }{\small \alpha }^{\prime }%
\end{array}%
\right) \wedge dx^{\prime \bar{A}}\right] .
\end{equation}%
Application of Stokes theorem yields%
\begin{equation}
\left[ dx^{A}{\LARGE (}_{\text{ }\wedge }^{\text{ }\cdot }{\Large )}dx^{i}%
\right] {\LARGE (}_{\text{ }\cdot }^{\text{ }\wedge }{\LARGE )}%
dx^{l}\int_{\partial E_{n}^{\prime }}\left( \frac{\partial }{\partial
x^{\prime l}}\frac{1}{r_{12}^{n-2}}\right) \left\{ dx^{\prime i}\cdot \left[
\left( 
\begin{array}{c}
{\small \delta }^{\prime }{\small \alpha }^{\prime } \\ 
{\small d}^{\prime }{\small \alpha }^{\prime }%
\end{array}%
\right) \wedge dx^{\prime \bar{A}}\right] \right\} ,  \label{50}
\end{equation}%
where we have indulged in the use of parentheses for greater clarity. It is
null if the differentiations of $\alpha $ go sufficiently fast to zero at
infinity.

The other term resulting from the integration by parts is 
\begin{equation}
-\left[ dx^{A}{\LARGE (}_{\text{ }\wedge }^{\text{ }\cdot }{\Large )}dx^{i}%
\right] {\LARGE (}_{\text{ }\cdot }^{\text{ }\wedge }{\LARGE )}%
dx^{l}\int_{E_{n}^{\prime }}\left( \frac{\partial }{\partial x^{\prime l}}%
\frac{1}{r_{12}^{n-2}}\right) \frac{\partial }{\partial x^{\prime i}}\left(
\left( 
\begin{array}{c}
{\small \delta }^{\prime }{\small \alpha }^{\prime } \\ 
{\small d}^{\prime }{\small \alpha }^{\prime }%
\end{array}%
\right) \wedge dx^{\prime \bar{A}}\right) .
\end{equation}%
This is zero because of cancellations that take place in groups of three
different indices. We shall devote the next subsection to dealing with the
intricacies of such cancellations.

In terms of Cartesian bases, we have, on the top line of the left hand side
of (42)%
\begin{equation}
dx^{A}\int_{E_{n}^{\prime }}\frac{(\delta ^{\prime }\alpha ^{\prime })\wedge
dx^{\prime \bar{A}}}{r_{12}^{n-2}}.
\end{equation}%
It is preceded by invariant operators, which we may ignore for present
purposes. We move $dx^{A}$ inside the integral, where we let $(\delta
^{\prime }\alpha ^{\prime })_{A}$ be the notation for the coefficients of $%
\delta ^{\prime }\alpha ^{\prime }.$ We thus have, for that first term,%
\begin{equation}
\int_{E_{n}^{\prime }}\frac{dx^{A}\wedge \lbrack (\delta ^{\prime }\alpha
^{\prime })_{A}dx^{\prime A}]\wedge dx^{\prime \bar{A}}}{r_{12}^{n-2}}.
\label{53}
\end{equation}%
The numerator can be further written as $(\delta ^{\prime }\alpha ^{\prime
})_{A}dx^{A}z^{\prime }.$ It is clear that $z$ and $(\delta ^{\prime }\alpha
^{\prime })_{A}dx^{\prime A}$ are invariants, but not immediately clear that 
$(\delta ^{\prime }\alpha ^{\prime })_{A}dx^{A}$ also is so$.$ Whether we
have the basis $dx^{A}$ or $dx^{\prime A}$ as a factor is immaterial. since
the invariance of $(\delta ^{\prime }\alpha ^{\prime })_{A}dx^{\prime A}$
can be seen as following from the matching of the transformations of $%
(\delta ^{\prime }\alpha ^{\prime })_{A}$ and $dx^{\prime A}$ each in
accordance with its type of covariance. The same matching applies if we
replace $\omega ^{\prime A}$ with $dx^{A}$, since $dx^{\prime A}$ and $%
dx^{A} $ transform in unison.

We have shown that (40) constitutes the decomposition of $\alpha $ into
closed and co-closed terms. Together with (41), it solves the problem of
integrating the system $d\alpha =\mu ,$ $\delta \alpha =\nu $, for given $%
\mu $ and $\nu $, and with the stated boundary condition

\subsection{Identical vanishing of some integrals}

As we are about to show, expressions (51) cancel identically (Expression
(50) cancels at infinity for fast vanishing, not identically).

Consider the first line in (51). Let $\alpha $ be of grade $h\geq 2$ (If $h$
were one, the dot product of $dx^{A}$ with $dx^{i}$ would be zero). Let $p$
and $q$ be a specific pair of indices in a given term in $\alpha $, i.e. in
its projection $a_{pqC}^{\prime },_{pq}dx^{A\text{ }}$upon some specific
basis element $dx^{A}$. Such a projection can be written as%
\begin{equation*}
(a_{pqC}^{\prime }dx^{\prime p}\wedge dx^{\prime q}\wedge dx^{\prime C},
\end{equation*}%
where $dx^{\prime A}$ is a unit monomial differential $1-$form (there is no
sum over repeated indices. We could also have chosen to write the same term
as%
\begin{equation*}
(a_{qpC}^{\prime }dx^{\prime q}\wedge dx^{\prime p}\wedge dx^{\prime C},
\end{equation*}%
with $a_{qpC}^{\prime }=-a_{pqC}^{\prime }.$ Clearly, $dx^{\prime C}$ is
uniquely determined if it is not to contain $dx^{p}$ and $dx^{q}.$ We then
have%
\begin{equation}
\delta ^{\prime }(a_{^{\prime }pqC}dx^{\prime p}\wedge dx^{\prime q}\wedge
dx^{\prime C})=a_{pqC}^{\prime },_{p}dx^{\prime q}\wedge dx^{\prime
C}-a_{pqC},_{q}dx^{\prime p}\wedge dx^{\prime C}.
\end{equation}%
The two terms on the right are two different differential $2-$forms. They
enter two different integrals, corresponding to $dx^{\prime q}\wedge
dx^{\prime C}$ and $dx^{\prime p}\wedge dx^{\prime C}$ components of $\delta
^{\prime }\alpha ^{\prime }$. To avoid confusion, we shall refer to the
basis elements in the integrals as $dx^{\prime B}$ since they are $(h-1)-$%
forms, unlike the $dx^{\prime A}$ of (54), which are differential $h-$forms

When taking the first term of (54) with $i=p$ into the top line of (51), the
factor at the front of the integral is%
\begin{equation*}
-(dx^{B}\cdot dx^{p})\wedge dx^{l}.
\end{equation*}%
But this factor is zero since $dx^{B}$ is $dx^{q}\wedge dx^{C}$, which does
not contain $dx^{P}$ as a factor. Hence, for the first term, we need only
consider $i=q.$ By the same argument, we need only consider $i=q$ for the
second term in (54). Upon multiplying the $dx^{\prime B}$'s by pertinent $%
dx^{\prime B}$'s, we shall obtain the combination%
\begin{equation*}
(a_{pqC}^{\prime },_{pq}-a_{pqC}^{\prime },_{qp})z^{\prime }
\end{equation*}%
as a factor inside the integral for the first line of (51). We could make
this statement because the factor outside also is the same one for both
terms: $(dx^{q}\wedge dx^{C})\cdot dx^{q}$ and $(dx^{p}\wedge dx^{C})\cdot
dx^{p}$ are equal. The contributions arising from the two terms on the right
hand side of (54) thus cancel each other out. We would proceed similarly
with any other pair of indices, among them those containing either $p$ or $q$%
. The annulment of the top line of (51) has been proved.

In order to prove the cancellation of the second line in (51), the following
considerations will be needed. A given $dx^{A}$ determines its corresponding 
$dx^{\prime \bar{A}}$, and vice versa. It follows then that only the term
proportional to $dx^{\prime A}$ in $d^{\prime }\alpha ^{\prime }$ exterior
multiplies $dx^{\prime \bar{A}}$, which is of the same grade as $d^{\prime
}\alpha ^{\prime }$, i.e. $h+1.$ Hence $dx^{A}\wedge dx^{i}$ is of grade $3$
or greater for $h>0.$ If $dx^{A}\wedge dx^{i}$ is not to be null, $dx^{i}$
cannot be in $dx^{A}$. Hence, $dx^{\prime \bar{A}}$ contains $dx^{i}$ as a
factor.

Let ($p,q,r$) be a triple of three different indices in $dx^{A}\wedge
dx^{i}. $ When $i$ is $p$ or $q$ or $r$, the respective pairs ($q,r$), ($r,p$%
) and ($p,q$) are in $dx^{A}$. We may thus write%
\begin{equation}
dx^{\prime A}=dx^{\prime C}\wedge dx^{\prime q}\wedge dx^{\prime r},\ \ \ \
\ \ dx^{\prime \bar{A}}=dx^{\prime p}\wedge dx^{\prime B}.
\end{equation}%
The coefficient of $dx^{\prime A}$ in $d^{\prime }\alpha ^{\prime }$ will be
the sum of three terms, one of which is%
\begin{equation}
(a_{Cr}^{\prime },_{q}-a_{Cq}^{\prime },_{r})dx^{\prime q}\wedge dx^{\prime
C}\wedge dx^{\prime r},  \label{56}
\end{equation}%
and the other two are cyclic permutations. We partial-differentiate (56)
with respect to $dx\prime ^{p}$ and multiply by $dx^{\prime p}\wedge
dx^{\prime B}$ on the right. We proceed similarly with $i=q$ and $i=r$, and
add all these contributions. We thus get%
\begin{equation}
(a_{Cr}^{\prime },_{qp}-a_{Cq}^{\prime },_{rp}+a_{Cp}^{\prime
},_{rq}-a_{Cr}^{\prime },_{pq}+a_{Cq}^{\prime },_{pr}-a_{Cp}^{\prime
},_{qr})z^{\prime }.
\end{equation}%
By virtue of equality of second partial derivatives, terms first, second and
third inside the parenthesis cancel with terms fourth, fifth and sixth. To
complete the proof, we follow the same process with another $dx^{\prime C}$
and the same triple ($p,q,r$) until we exhaust all the options. We then
proceed to choose another triple and repeat the same process until we are
done with all the terms, which completes the proof of identical vanishing of
the second term arising from one of the two integrations by parts of the
previous subsection.

\section{Hodge's Theorems}

The ``beyond'' in the title of this section responds to the fact that we
shall be doing much more than reproducing Hodge's theorem. As is the case
with Helmholtz theorem, we are able to specify in terms of integrals what
the different terms are.

We shall later embed Riemannian spaces $R_{n}$ in Euclidean spaces $E_{N}$,
thus becoming $n-$ surfaces. As an intermediate step, we shall apply the
traditional Helmholtz approach to regions of Euclidean spaces, i.e. $R_{n}$%
's ab initio embedded in $E_{n}$. The harmonic form ---which is of the
essence in Hodge's theorem--- emerges from the Helmholtz process in the new
venues.

\subsection{Transition from Helmholtz to Hodge}

Though visualization is not essential to follow the argument, it helps for
staying focused. For that reason, we shall argue in 3-D Euclidean space. It
does not interfere with the nature of the argument.

On a region $R$ of $E_{3}$, including the boundary, define a differential $%
1- $form or $2-$form $\alpha $. Let $A$ denote any continuously
differentiable prolongation of $\alpha $ that vanish sufficiently fast at
infinity. On $R$, we have $dA=d\alpha $ and $\delta A=\delta \alpha $. We
can apply Helmholtz theorem to the differential forms $A$. In order to
minimize clutter, we write it in the form%
\begin{eqnarray}
-4\pi A &=&d...\int_{R^{\prime }}\frac{\delta ^{\prime }A^{\prime }...}{%
r_{12}}+\delta \int_{R^{\prime }}...\frac{d^{\prime }A^{\prime }...}{r_{12}}+
\notag  \label{58} \\
&&+d...\int_{E_{3}^{\prime }-R^{\prime }}\frac{\delta ^{\prime }A^{\prime
}...}{r_{12}}+\delta ...\int_{E_{R}^{\prime }-R^{\prime }}\frac{d^{\prime
}A^{\prime }...}{r_{12}},
\end{eqnarray}%
where $r_{12}=[(x-x^{\prime })^{2}+(y-y^{\prime })^{2}+(z-z^{\prime
})^{2}]^{1/2}.$ We shall keep track of the fact, at this point obvious, that
in the first two integrals on the right, $r^{\prime }$ is in $R^{\prime }$.
It is outside $R^{\prime }$ in the other two integrals, which will depend on
the prolongation. By representing those terms simply as $\mathcal{F}$, we
have%
\begin{equation}
-4\pi A=d...\int_{R^{\prime }}\frac{\delta ^{\prime }\alpha ^{\prime }...}{r}%
+\delta ...\int_{R^{\prime }}\frac{d^{\prime }\alpha ^{\prime }...}{r}+%
\mathcal{F}.  \label{59}
\end{equation}%
Since these equations yield $A$ everywhere in $E_{3}$ (i.e. $r$ not limited
to $R$), they yield in particular what $A$ and $\mathcal{F}$ are in $R$ . We
can thus write%
\begin{equation}
-4\pi \alpha =d...\int_{R^{\prime }}\frac{\delta ^{\prime }\alpha ^{\prime
}...}{r}+\delta ...\int_{R^{\prime }}\frac{d^{\prime }\alpha ^{\prime }...}{r%
}+\mathcal{F},
\end{equation}%
$\mathcal{F}$ not having changed except that $\mathcal{F}$ in (60) refers
only to what it is in $R$ but it remains a sum of integrals in $%
E_{3}^{\prime }-R.$ The prolongations will be determined as different
solutions of a differential system to be obtained as follows.

By following the same process as in Helmholtz theorem, we obtain, in
particular, 
\begin{equation}
-4\pi d\alpha =d\delta ...\int_{R^{\prime }}\frac{d^{\prime }\alpha ^{\prime
}...}{r}+d\mathcal{F},  \label{61}
\end{equation}%
and similarly for $-4\pi \delta \alpha $ (just exchange $d$ and $\delta $).

Now, the first term on the right hand side of (61) will not become simply $%
-4\pi d\alpha $ as was the case in the previous section. It will yield two
terms. One of them is $-4\pi d\alpha $, and the other one is made to cancel
with $d\mathcal{F}$, thus determining a differential equation to be
satisfied by $\mathcal{F}$. To this we have to add another differential
equation arising from application of $\delta $ to (60). Together they
determine the differential system to be determined by $\mathcal{F}.$ Thus $%
-4\pi \alpha $ will be given by the three term decomposition (60). Notice
that, in the process, we avoid integrating over $E_{3}^{\prime }-R^{\prime }$
and instead solving a differential system in $R$, since the left hand side
and the first term on the right hand side of (61) pertain to $\alpha $.

From now one, we shall make part of the theorems that the prolongations are
solutions of a certain differential systems, later to be made explicit.

\subsection{Hodge theorem in regions of E$_{n}$}

Let $\alpha $ be a differential $k-$form satisfying the equations $d\alpha
=\mu $ and $\delta \alpha =\nu $, and given at the boundary of a region of $%
E_{n}$. We proceed to integrate this system. (60) now reads 
\begin{equation}
-(n-2)S_{n-1}\alpha =d\left[ \int_{R^{\prime }}\frac{\delta ^{\prime }\alpha
^{\prime }...}{r_{12}}\right] +\delta \left[ \int_{R^{\prime }}\frac{%
d^{\prime }\alpha ^{\prime }...}{r_{12}}\right] +\mathcal{F},
\end{equation}%
where $R$ is a region of Euclidean space that contains the origin and where $%
r_{12}$ is the magnitude of the Euclidean distance between hypothetical
points of components ($x,y,...u,v$).and ($x^{\prime },y^{\prime
},...u^{\prime },v^{\prime }$), all the coordinates chosen as Cartesian to
simplify visualization. We said hypothetical because the interpretation as
distance only makes sense when we superimpose $E_{n}$ and $E_{n}^{\prime }$.

When we apply either $d$ or $\delta $ to (62), we shall use, as before, $%
d\delta +$ $\delta d=\partial \partial ,$ with one of the terms on the left
moved to the right ($d\delta =...,$ $\delta d=...$ respectively). By
developing the $\partial \partial $ term, it becomes the same as term on the
right (i.e $d\alpha $ or $\delta \alpha $). It will cancel with the term on
the left. The terms that vanished identically also vanish now, precisely
because this is an identical vanishing. We are thus left with the total
differential terms. If apply Stokes theorem, as before. these terms no
longer disappear at the boundary. Hence, we are left with the two equations%
\begin{equation}
\left[ dx^{A}{\LARGE (}_{\text{ }\wedge }^{\text{ }\cdot }{\Large )}dx^{i}%
\right] {\LARGE (}_{\text{ }\cdot }^{\text{ }\wedge }{\LARGE )}%
dx^{l}\int_{R^{\prime }}\left( \frac{\partial \frac{1}{r_{12}^{n-2}}}{%
\partial x^{\prime l}}\right) dx^{\prime i}\cdot \left[ \left( 
\begin{array}{c}
{\small \delta }^{\prime }{\small \alpha }^{\prime } \\ 
{\small d}^{\prime }{\small \alpha }^{\prime }%
\end{array}%
\right) \wedge dx^{\prime \bar{A}}\right] +\left( 
\begin{array}{c}
d \\ 
\delta%
\end{array}%
\right) \mathcal{F}=0  \label{63}
\end{equation}%
(Refer to (50)). Hence, the solution to Helmholtz problem is given by the
pair of equations (62)-(63).

We shall now show that $\mathcal{F}$ is harmonic, i.e. $(d\delta +\delta d)%
\mathcal{F}=0.$.We shall apply $\delta $ and $d$ to the first and second
lines of (63). Start by rewriting the first terms in (63) in the form, (49),
they took before applying Stoke's theorem. Upon applying the $\delta $
operator to the first line, we have, for $\delta d\mathcal{F}$,%
\begin{equation}
dx^{h}\cdot \lbrack (dx^{A}\cdot dx^{i})\wedge dx^{l}]\int_{R^{\prime }}%
\frac{\partial ^{2}}{\partial x^{\prime h}\partial x^{\prime i}}\left[
\left( \frac{\partial }{\partial x^{\prime l}}\frac{1}{r_{12}^{n-2}}\right)
\left( 
\begin{array}{c}
{\small \delta }^{\prime }{\small \alpha }^{\prime } \\ 
{\small d}^{\prime }{\small \alpha }^{\prime }%
\end{array}%
\right) \wedge dx^{\prime \bar{A}}\right] .
\end{equation}%
Since this term happens to vanish, the computation will take place up to the
factor $-1$, provided it is common to all terms in a development into
explicit terms. We do so because (64) will be shown to vanish identically.

For $dx^{h}\cdot \lbrack (dx^{A}\cdot dx^{i})\wedge dx^{l}]$ to be different
from zero, $h$ and $i$ must be different and contained in $A.$ Since $dx^{l}$
is not in $dx^{A}$, the product $dx^{h}\cdot dx^{l}$ is zero. Hence%
\begin{equation}
dx^{h}\cdot \lbrack (dx^{A}\cdot dx^{i})\wedge dx^{l}]=[dx^{h}\cdot
(dx^{A}\cdot dx^{i})]\wedge dx^{l}].
\end{equation}%
We can always write $dx^{A}$ as%
\begin{equation}
dx^{h}\wedge dx^{j}\wedge dx^{C}\wedge dx^{i}.
\end{equation}%
This is antisymmetric in the pair ($i,h$), which combines with the symmetry
inside the integral to annul this term. Notice that we did not have to
assign specific values for ($i,h$), but we had to ``go inside'' $dx^{A}$. We
mention this for contrast with the contents for the next paragraph. We have
proved so far that $\delta d\mathcal{F}=0.$

We proceed to prove that $d\delta \mathcal{F}=0$. We rewrite the left hand
side of (63) as in (49) and proceed to apply $d$ to it. We shall now have%
\begin{equation}
dx^{h}\wedge \lbrack (dx^{A}\wedge dx^{i})\cdot dx^{l}]\int_{R^{\prime }}%
\frac{\partial ^{2}}{\partial x^{\prime h}\partial x^{\prime i}}\left[
\left( \frac{\partial }{\partial x^{\prime l}}\frac{1}{r_{12}^{n-2}}\right) 
{\small d}^{\prime }{\small \alpha }^{\prime }\wedge dx^{\prime \bar{A}}%
\right] .  \label{67}
\end{equation}%
It is clear that, when $l$ takes a value different from the value taken by $%
i $, we again have cancellation due to the same combination of
antisymmetry-symmetry as before. But the terms $dx^{i}\cdot dx^{l}$ would
seem to interfere with the argument, but it does not. We simply have to be
more specific than before with the groups of terms that we put together. We
put together only terms where we have $dx^{r}\wedge dx^{s}$ arising from ($%
h=r,$ $i=s$) and ($h=s,$ $i=r$). When the running index $l$ takes the values 
$r$ or $s$, the resulting factor at the front of the integral will belong to
a different group. We have thus shown that (67) cancels out and, therefore, $%
d\delta \mathcal{F}=0.$ To be precise, we have not only proved that $%
\mathcal{F}$ is harmonic, but that it is ``hyper-harmonic'', meaning
precisely that: $\delta d\mathcal{F}=0$ and $d\delta \mathcal{F}=0.$

\subsection{Hodge's theorem in hypersurfaces of E$_{N}$}

A manifold embedded in a Euclidean space of the same dimension will be
called a region thereof. A hypersurface is a manifold of dimension $n$
embedded in a Euclidean space $E_{N}$ where $N>n$. The treatment here is the
same as in subsection 6.1, the hypersurface playing the role of the region.
The only issue that we need to deal with is a practical one having to do
with the experience of readers. Helmholtz magnificent theorem belongs to an
epoch where vector (and tensor) fields often took the place of differential
forms. This can prompt false ideas as we now explain.

Let $\mathbf{v}$ be a vector field $\mathbf{v\equiv }$ $a^{\lambda }(u,v)%
\mathbf{\hat{a}}_{\lambda }$ ($\lambda =1,2$) on a surface $x^{i}(u,v)$ ($%
i=1,2,3$) embedded in $E_{3}$, the frame field $\mathbf{\hat{a}}_{\lambda }$
being orthonormal. It can be tangent or not tangent. By default, the vector
field is zero over the remainder of $E_{3}$. In its present form , Helmholtz
theorem would not work for this field since the volume integrals over $E_{3}$
would be zero. This is a spurious implication because the theorem should be
about algebras of differential forms, not tangent spaces.

Let $\mu $ be the differential $1-$form$\ a_{\lambda }(u,v)\hat{\omega}%
^{\lambda }$, the basis $\hat{\omega}^{\lambda }$ being dual to the constant
orthonormal basis field $\mathbf{a}_{\lambda }$. This duality yields $%
a_{\lambda }=a^{\lambda }$. No specific curve is involved in the definition
of $\mu $, which is a function of curves, function determined by its
coefficients $a_{\lambda }(u,v)$ \ The specification of a vector field on a
surface, $\mathbf{v}$, on the other hand needs to make reference to a
surface for its definition. And yet the components of $d\mu $ and $\delta
\mu $ (which respectively are a $2-$form and a $0-$form) enter non-null
volume integrals, which pertain to $3-$forms . The fact that most components
(in the algebra) of an $k-$form are zero is totally irrelevant. The
Helmholtz theorem for, say, a differential $1-$form $\mu $ can be formulated
in any sufficiently high dimensional Euclidean space regardless of whether
the ``associated'' vector field $\mathbf{v}$ is zero outside some surface.

Similarly, Helmholtz theorem for a differential $n-$form in $E_{N}$ involves
the integration of differential $N-$forms, built upon the interior
differential $(n-1)-$form and the exterior differential $(n+1)-$form. In
considering simple examples (say a plane in $3-$space), one can be misled or
confused if one does not take into account the role of $1/r$, or else we
might be obtaining an indefinite integral. Assume finite $\int \lambda
(x,y)dx\wedge dy$ when integrating over the $xy$ plane. The integral $\int
\lambda (x,y)dx\wedge dy\wedge dz$ would be divergent, but need not be so if
there is some factor that goes to zero sufficiently fast at infinity of $z$
and $-z.$

It is also worth mentioning that ---in the case of a hypersurface like in
the case of a surface in $E_{3}$---, $r_{12}$ represents a chord, which is
not in the hypersurface.

A final issue worth addressing is the following. If the Laplacian of the
appropriate power of $1/r$ is zero and it multiplies the differential form
outside the region or outside the hypersurface, why does the prolongation
make a difference. This is a pseudo-problem easy to understand already at
the point of equations (2) and (3), i.e. before we deal with prolongations.
Those Laplacians are not generalizations of functions and cannot, therefore,
be treated as such \cite{Schwarz}.

We conclude this subsection with the observation that, in Helmholtz theorems
for differential forms, the variety of disconnected concepts that enter
Helmholtz theorem for vector fields on surfaces in $E_{3}$ (a vector field,
a surface, a gradient, a divergence, a curl, integrands and $E_{3}$) merge
into or directly connect with the concept of differential form as an
integrand in the Helmholtz theorem of the K\"{a}hler calculus, where we have

(a) a \textbf{differential }$1-$\textbf{form}, in lieu of a \textit{vector
field},

(b) from which we construct through K\"{a}hler differentiation \textbf{an
inhomogeneous differential form, }in lieu of \textit{divergence, gradient
and curl},

(c) from which in turn we build and evaluate (read integrate) a \textbf{%
differential 3-form}, in lieu of three \textit{volume integrals, one each
for the components of the vector field},

(d) and we restrict the \textbf{coefficients of the differential forms }in
(a) and (b)\textbf{\ }to surfaces, \textit{in lieu of vector fields defined
on surfaces}.

\subsection{Hodge theorems in Riemannian spaces}

We shall consider a Helmholtz-Hodge extension of Hodge's theorem (i.e. a
theorem of integration) and the standard Hodge theorem, which is a
consequence of the former.

Consider now a differentiable manifold $R_{n}$ endowed with a Euclidean
metric. By the Schl\"{a}fli-Janet-Cartan theorem \cite{Schlafli},\cite{Janet}%
,\cite{Cartan27}, it can be embedded in a Euclidean space of dimension $%
N=n(n+1)/2.$ Hence, a Helmholtz-Hodge theorem follows for orientable
Riemannian manifolds that satisfy the conditions for application of Stokes
theorem by viewing them as hypersurfaces in Euclidean spaces. At this point
in our argument, the positive definiteness of the metric is required, or
else we would have to find a replacement for the Laplacians considered in
previous sections. The result is local, meaning non global, remark made in
case the term local might send some physicists in a different direction. For
clarity, the evaluation of the Laplacian now satisfies 
\begin{equation}
1=\frac{1}{(N-2)S_{N-1}}\int_{E_{N}}\partial \partial \frac{1}{r^{N-2}}z,
\end{equation}%
where $r$ is the radial coordinate in $N-$dimensional space. Needless to say
that it also applies to regions and hypersurfaces of $E_{N}$ that contain
the origin. As a consequence of the results in the previous subsections, we
have the following.

\textbf{Helmholtz-Hodge theorem:}

Hodge's theorem is constituted by Eqs. (69)-(71): For differential $k-$forms
in Riemannian spaces $R_{n}$%
\begin{equation}
-(N-2)S_{N-1}\alpha =d\left[ \omega ^{A}\int_{R_{n}^{\prime }}\frac{(\delta
^{\prime }\alpha ^{\prime })\wedge \omega ^{\prime \bar{A}}}{r_{12}^{N-2}}%
\right] +\delta \left[ \omega ^{A}\int_{R_{n}^{\prime }}\frac{(d^{\prime
}\alpha ^{\prime })\wedge \omega ^{\prime \bar{A}}}{r_{12}^{N-2}}\right] +%
\mathcal{F},
\end{equation}

\begin{equation}
\left( 
\begin{array}{c}
d \\ 
\delta%
\end{array}%
\right) \mathcal{F=-}\left[ dx^{A}{\LARGE (}_{\text{ }\wedge }^{\text{ }%
\cdot }{\Large )}dx^{i}\right] {\LARGE (}_{\text{ }\cdot }^{\text{ }\wedge }%
{\LARGE )}dx^{l}\int_{R_{n}^{\prime }}\left( \frac{\partial \frac{1}{%
r_{12}^{n-2}}}{\partial x^{\prime l}}\right) dx^{\prime i}\cdot \left[
\left( 
\begin{array}{c}
{\small \delta }^{\prime }{\small \alpha }^{\prime } \\ 
{\small d}^{\prime }{\small \alpha }^{\prime }%
\end{array}%
\right) \wedge dx^{\prime \bar{A}}\right] ,  \label{70}
\end{equation}%
with $r_{12}$ being defined in any Euclidean space of dimension $N\geq
n(n+1)/2$ where we consider $R_{n}$ to be embedded.

As previously discussed, $r_{12}$ represents a chord. We insist once more
that $\omega ^{\prime \bar{A}}$ is determined by the specific term in $%
\delta ^{\prime }\alpha ^{\prime }$ and $d^{\prime }\alpha ^{\prime }$ that
it multiplies. $\mathcal{F}$ is undetermined by solutions of the system $%
\delta \alpha =0$, $d\alpha =0.$ So is, therefore $\alpha .$

Hodge's theorem, as opposed to Helmholtz-Hodge theorem, is about
decomposition. Hence, once again, uniqueness refers to something different
from the uniqueness in the theorem of subsection (3.2), which refers to a
differential system.

One might be momentarily tempted to now apply (69) to (70). We would get an
identity, $\mathcal{F=F}$, by virtue of the orthogonality of the subspace of
the harmonic differential forms to the subspaces of closed and co-closed
differential forms.

\textbf{Hodge's theorem:}

Any differential $k-$form, whether of homogeneous grade or not, can be
uniquely decomposed into closed, co-closed and hyper-harmonic terms. For
differential $k-$forms, the theorem is an immediate consequence of (69). For
differential forms which are not of homogeneous grade, the theorem also
applies because one only needs to add the decompositions of the theorem for
the different homogeneous $k-$forms that constitute the inhomogeneous
differential form.

\section{Concluding Remarks.}

We have obtained by computation results that, to our knowledge, had not been
addressed before. We have gone as far as obtaining a Helmholtz theorem for
differential $k-$forms in Riemannian manifolds. By virtue of the acquisition
in general of a third, harmonic term, we have come to call it
Helmholtz-Hodge, which goes far beyond Hodge's theorem.

The more general results may still keep us far away from practical
applications of these theorems, except in isolated cases in low dimensions.
The reason is that the solution of problems of embedding are not trivial. It
is a symptom of their difficulty that the embedding of a 3-D manifold in a
6-D Euclidean space takes the last pages of Cartan's treatise of integration
of exterior systems, a book dedicated to the Cartan-K\"{a}hler theory \cite%
{Cartan34}.

Of more practical interest is the fact that these results show the
tremendous potential of the K\"{a}hler calculus, both in physics and
mathematics. It would be too self-serving to mention here specific results.
Interested readers can go to the authors web site
www.physical-unification.com for references.

\section{Acknowledgements}

Conversations with Prof. Z. Oziewicz are acknowledged. Funding from PST\
Associates is deeply appreciated.

\end{document}